\documentclass[11pt,a4paper,reqno]{amsart}

\usepackage[english]{babel}
\usepackage[latin1]{inputenc}
\usepackage[hidelinks]{hyperref}
\usepackage{amsmath,
math rsfs,
xcolor,
dsfont
}

\usepackage{geometry}
\newtheorem{theorem}{Theorem}[section]

\newtheorem{conjecture}{Conjecture}[section]
\theoremstyle{remark}
\newtheorem{remark}{Remark}[section]

\begin{document}
 
\title[On the extensions of the De Giorgi approach]{On the extensions of the De Giorgi approach \\ to nonlinear hyperbolic equations}

\author[L. Tentarelli]{Lorenzo Tentarelli$^1$}
\address{Sapienza Universit\`a di Roma, Dipartimento di Matematica, Piazzale Aldo Moro, 5, 00185, Roma, Italy.}
\email{tentarelli@mat.uniroma1.it}

\footnotetext[1]{Author supported by the FIR grant 2013 ``Condensed Matter in Mathematical Physics (Cond-Math)'' (code RBFR 13WAET).}

\begin{abstract} 
 In this talk we present an overview on the extensions of the De Giorgi approach to general second order nonlinear hyperbolic equations. We start with an introduction to the original conjecture by E. De Giorgi (\cite{degiorgi:1,degiorgi:2}) and to its solution by E. Serra and P. Tilli (\cite{serra:2}). Then, we discuss a first extension of this idea (Serra$\&$Tilli, \cite{serra:3}) aimed at investigating a wide class of homogeneous equations. Finally, we announce a further extension to nonhomogeneous equations, obtained by the author in \cite{tentarelli} in collaboration with P. Tilli. 
\end{abstract}

\maketitle

\tableofcontents

\section{De Giorgi's conjecture.}

In 1996, E. De Giorgi stated the following conjecture on weak solutions of the {\it defocusing} NLW equation.

\begin{conjecture}[De Giorgi, \cite{degiorgi:1,degiorgi:2}]
 \label{conj:degiorgi}
 Let $w_0,w_1\in C_0^{\infty}(\mathbb{R}^n)$, let $k>1$ be an integer; for every positive real number $\varepsilon$, let $w_{\varepsilon}=w_{\varepsilon}(t,x)$ be the minimizer of the functional
 \begin{equation}
  \label{eq:functional_degiorgi}
  F_\varepsilon(u):=\int_0^{\infty}\int_{\mathbb{R}^n}e^{-t/\varepsilon}\,\left(|u''(t,x)|^2+\tfrac{1}{\varepsilon^2}|\nabla u(t,x)|^2+\tfrac{1}{\varepsilon^2}|u(t,x)|^{2k}\right)\,dx\,dt
 \end{equation}
 in the class of all $u$ satisfying the initial conditions
 \begin{equation}
  \label{eq:initial_conditions}
  u(0,x)=w_0(x),\qquad u'(0,x)=w_1(x).
 \end{equation}
 Then, there exists $\displaystyle\lim_{\varepsilon\downarrow0}w_{\varepsilon}(t,x)=w(t,x)$, satisfying the equation
 \begin{equation}
  \label{eq:wave_degiorgi}
  w''=\Delta w-kw^{2k-1}.
 \end{equation}
\end{conjecture}

\begin{remark}
 In the statement of the conjecture, we maintained the original formulation of \cite{degiorgi:1,degiorgi:2} and we only changed notation, according to that we use in the sequel. The same thing holds for all the results we mention in this paper. In addition, we recall that $u'(t,x)$ denotes $\frac{\partial u}{\partial t}(t,x)$ and that, for the sake of simplicity, we always omit the dependence of the functional spaces on $\mathbb{R}^n$, i.e. $H^1=H^1(\mathbb{R}^n)$, $L^p=L^p(\mathbb{R}^n)$ and so on.
\end{remark}

In order to better understand the meaning of the conjecture, it is worth stressing some characteristic features of the functional $F_\varepsilon$.

First we note that it involves second order time derivatives. Thus, a minimizer of $F_\varepsilon$ solves a fourth order PDE. However, if one computes the formal {\it Euler--Lagrange} equation satisfied by a minimizer $w_{\varepsilon}$, then one obtains
\[
 \varepsilon^2(e^{-t/\varepsilon}\,w_{\varepsilon}'')''=e^{-t/\varepsilon}\,(\Delta w_{\varepsilon}-kw_{\varepsilon}^{2k-1})
\]
and thus, expanding and dropping $e^{-t/\varepsilon}$,
\begin{equation}
 \label{eq:euler_lagrange}
 \varepsilon^2\,w_{\varepsilon}''''-2\varepsilon w_{\varepsilon}'''+w_{\varepsilon}''=\Delta w_{\varepsilon}-kw_{\varepsilon}^{2k-1}.
\end{equation}
Consequently, if one assumes that $w_{\varepsilon}\to w$ in some suitable sense and lets $\varepsilon\downarrow0$, then one formally obtains \eqref{eq:wave_degiorgi}.

On the other hand, we also remark that, as $F_{\varepsilon}$ is defined through integrals over the ``space--time'' $[0,\infty)\times\mathbb{R}^n$, the initial conditions of the Cauchy problem are in fact {\it boundary conditions} for the minimization problem.

In addition, it is convenient to stress the singular nature of the {\it integration weight} $e^{-t/\varepsilon}$. More precisely, one can see that $\varepsilon^{-1}\,e^{-t/\varepsilon}\,dt$ is an approximate {\it Dirac delta} measure and hence, at least formally,
\[
 \varepsilon F_{\varepsilon}(u)\approx\int_{\mathbb{R}^n}\big(|\nabla w_0(x)|^2+|w_0(x)|^{2k}\big)\,dx,\qquad\text{as}\quad\varepsilon\downarrow0.
\]
Hence, this prevents a straightforward application of classical techniques of variational convergence, such as $\Gamma$--{\it convergence}. The previous asymptotic expansion, indeed, shows that this technique does not provide useful information on the limit behavior of the sequence of the minimizers.

Finally, one can note that $F_\varepsilon$ is convex (for fixed $\varepsilon>0$) and that therefore, up to some suitable technical adaptation, the proof of the existence and uniqueness of the minimizers is not a demanding issue.

\medskip
We also recall that the existence of \emph{global} solutions for the Cauchy problem \eqref{eq:wave_degiorgi}$\&$\eqref{eq:initial_conditions} is not new (see e.g. \cite{strauss} and the references therein). Actually, as highlighted in \cite{serra:1,serra:2}, the originality of the strategy hinted by De Giorgi lies in \emph{how} he intended to exploit techniques from the Calculus of Variations. The variational approaches to the wave equation $w''=\Delta w$ and its nonlinear variants that can one can find in the literature (see e.g. \cite{strauss,struwe} and references therein) are based on the interpretation of $w''=\Delta w$ as the Euler--Lagrange equation of the functional
\[
 I(w):=\int_0^{\infty}\int_{\mathbb{R}^n}\left(|w'(t,x)|^2-|\nabla w(t,x)|^2\right)\,dx\,dt
\]
(with possibly lower order terms like $|w|^{2k}$). However, since $I$ is neither convex nor bounded from below, one is forced to search for \emph{critical points} rather than \emph{global minimizers}. Unfortunately, functionals like $I$ behave badly also for the application of Critical Point Theory, so that only partial results can be proved. De Giorgi, on the contrary, introduces a new functional $F_\varepsilon$ that is quite easy to minimize (regardless of the magnitude of $k$) and thus moves the problem to the investigation of the limit behavior of the sequence of the minimizers.

\section{The proof of the conjecture.}

In 2012, E. Serra and P. Tilli showed that Conjecture \ref{conj:degiorgi} is in fact true. Precisely, in \cite{serra:2}, they proved the following theorem.

\begin{theorem}[Serra$\&$Tilli, \cite{serra:2}]
 \label{theorem:serra-tilli-uno}
 For $p\geq2$ and $\varepsilon>0$, let $w_{\varepsilon}(t,x)$ denote the unique minimizer of the strictly convex functional
 \[
  F_\varepsilon(u)=\int_0^{\infty}\int_{\mathbb{R}^n}e^{-t/\varepsilon}\,\left(|u''(t,x)|^2+\tfrac{1}{\varepsilon^2}|\nabla u(t,x)|^2+\tfrac{1}{\varepsilon^2}|u(t,x)|^{p}\right)\,dx\,dt
 \]
 under the boundary conditions \eqref{eq:initial_conditions}, where $w_0$ and $w_1$ are given functions such that
 \[
  w_0,\,w_1\in H^1\cap L^p.
 \]
 Then:
 \begin{itemize}
  \item[(a)] Estimates. There exists a constant $C$ (which depends only on $w_0,\,w_1,\,p$ and $n$) such that, for every $\varepsilon\in(0,1)$,
   \begin{gather*}
    \displaystyle\int_0^T\int_{\mathbb{R}^n}\left(|\nabla w_{\varepsilon}(t,x)|^2+|w_{\varepsilon}(t,x)|^p\right)\,dx\,dt\leq CT,\qquad\forall T>\varepsilon,\\[.3cm]
    \displaystyle\int_{\mathbb{R}^n}|w_{\varepsilon}'(t,x)|^2\,dx\leq C\quad\mbox{and}\quad\int_{\mathbb{R}^n}|w_{\varepsilon}(t,x)|^2\,dx\leq C(1+t^2\,),\qquad\forall t\geq0,
   \end{gather*}
  and, for every function $h\in H^1\cap L^p$
  \[
   \left|\int_{\mathbb{R}^n}w_{\varepsilon}''(t,x)h(x)\,dx\right|\leq C\left(\|h\|_{L^p}+\|\nabla h\|_{L^2}\right),\qquad\mbox{for a.e.}\quad t>0.
  \]
  \item[(b)] Convergence. Every sequence $w_{\varepsilon_i}$ (with $\varepsilon_i\downarrow0$) admits a subsequence which is convergent, in the strong topology of $L^q((0,T)\times A)$ for every $T>0$ and every bounded open set $A\subset\mathbb{R}^n$ (with arbitrary $q\in[2,p)$ if $p>2$ and $q=p$ if $p=2$), almost everywhere in $\mathbb{R}^+\times\mathbb{R}^n$ and in the weak topology of $H^1((0,T)\times\mathbb{R}^n)$ for every $T>0$, to a function $w$ such that
  \[
   \begin{array}{ll}
    w\in L^{\infty}(\mathbb{R}^+;L^p), & \nabla w\in L^{\infty}(\mathbb{R}^+;L^2),\\[.5cm]
    w'\in L^{\infty}(\mathbb{R}^+;L^2), & w\in L^{\infty}((0,T);H^1)\quad\forall T>0,
   \end{array}
  \]
  which solves in $\mathbb{R}^+\times\mathbb{R}^n$ the nonlinear wave equation 
  \begin{equation}
   \label{eq:wave_serra-tilli}
   w''=\Delta w-\tfrac{p}{2}|w|^{p-2}\,w
  \end{equation}
  with initial conditions as in \eqref{eq:initial_conditions}.
  \item[(c)] Energy inequality. Letting
  \[
   \mathcal{E}(t):=\int_{\mathbb{R}^n}\left(|w'(t,x)|^2+|\nabla w(t,x)|^2+|w(t,x)|^p\right)\,dx,
  \]
  the function $w(t,x)$ satisfies the energy inequality
  \[
   \mathcal{E}(t)\leq\mathcal{E}(0)=\int_{\mathbb{R}^n}\left(|w_1(x)|^2+|\nabla w_0(x)|^2+|w_0(x)|^p\right)\,dx,\qquad\mbox{for a.e.}\quad t>0.
  \]
 \end{itemize}
\end{theorem}

\begin{remark}
 We stress the fact that, in \emph{(b)}, the limit function $w$ solves \eqref{eq:wave_serra-tilli} in a \emph{distributional} (or \emph{weak}) sense, namely
 \begin{align*}
  \int_0^{\infty}\int_{\mathbb{R}^n}w'(t,x)\varphi'(t,x)\,dx\,dt= & \, \int_0^{\infty}\int_{\mathbb{R}^n}\nabla w(t,x)\cdot\nabla\varphi(t,x)\,dx\,dt+\\[.2cm]
                                                                  & \, +\int_0^{\infty}\int_{\mathbb{R}^n}\tfrac{p}{2}|w(t,x)|^{p-2}\,w(t,x)\varphi(t,x)\,dx\,dt
 \end{align*}
 for every $\varphi\in C_0^{\infty}(\mathbb{R}^+\times\mathbb{R}^n)$. In the sequel we only deal with this type of solutions.
\end{remark}

Some comments are in order. First, \textsc{Conjecture} \ref{conj:degiorgi} deals with the nonlinearity $|w|^{2k}$ with $k$ \emph{integer}, while \textsc{Theorem} \ref{theorem:serra-tilli-uno} treats $|w|^{p}$ without the assumption of $p$ integer. Another relevant feature of \textsc{Theorem} \ref{theorem:serra-tilli-uno}, is that the assumptions on the initial data $w_0,\,w_1$ are much weaker than those of the conjecture.

On the other hand, the convergence of the sequence of the minimizers is obtained up to extracting subsequences, thus ``losing'' the \emph{uniqueness} claimed in the conjecture. In particular, it is an open problem to avoid the extraction of subsequences when $p$ is large.

In addition, \textsc{Theorem} \ref{theorem:serra-tilli-uno} establishes an estimate for the \emph{mechanical energy} $\mathcal{E}$ usually associated with \eqref{eq:wave_serra-tilli}, which proves that the obtained solutions are \emph{of energy class} in the sense of Struwe (see \cite{struwe}). When $p$ is ``sufficiently'' small, the inequality is in fact an equality, whereas, when $p$ is large, energy \emph{conservation} is still open.

For the sake of completeness, we mention that \cite{stefanelli} discusses a simplified version of the conjecture on bounded intervals. However, that paper only deals with the proof of \eqref{eq:wave_serra-tilli} and does not treat the fulfillment of the initial condition $w'(0,x)=w_1(x)$.

\section{Extension to homogeneous equations.}

Now, one can easily see that, setting
\[
 \mathcal{W}(v)=\int_{\mathbb{R}^n}\bigg(\frac{1}{2}|\nabla v|^2+\frac{1}{p}|v|^p\bigg)\,dx,
\]
up to some multiplicative constants equation \eqref{eq:wave_serra-tilli} reads
\begin{equation}
 \label{eq:wave_general}
 w''(t,x)=-\nabla\mathcal{W}(w(t,\cdot))(x),
\end{equation}
where $\nabla\mathcal{W}$ denotes the {\it G\^ateaux derivative} of the functional $\mathcal{W}$. Therefore, it is natural to wonder if the sequence of the minimizers of the functional $F_\varepsilon$, that here is defined by
\begin{equation}
 \label{eq:functional_homogeneous}
 F_\varepsilon(u):=\int_0^te^{-t/\varepsilon}\bigg(\int_{\mathbb{R}^n}\frac{\varepsilon^2|u''(t,x)|}{2}\,dx+\mathcal{W}(u(t,\cdot))\bigg)\,dt,
\end{equation}
converges to a solution of the Cauchy problem associated with \eqref{eq:wave_general}, even for different choices of $\mathcal{W}$.

\begin{remark}
 In \eqref{eq:functional_homogeneous} one uses a different scaling in $\varepsilon$, with respect to \eqref{eq:functional_degiorgi}. This is due to the fact that in the {\it abstract} framework this choice simplifies computations. However, this does not yield significant differences.
\end{remark}

This problem has been solved again by E. Serra and P. Tilli, in \cite{serra:3}. Before showing the statements of the main results, it is necessary to point out under which assumptions on the functional $\mathcal{W}$ (that we refer to as assumption {\bf (H)} in the following), they are valid.
\begin{itemize}
 \item[{\bf (H)}] The functional $\mathcal{W}:L^2\to[0,\infty]$ is lower semi--continuous in the weak topology of $L^2$, i.e
 \[
  \mathcal{W}(v)\leq\liminf_{k}\mathcal{W}(v_k),\qquad\mbox{whenever}\quad v_k\rightharpoonup v\quad\mbox{in}\quad L^2.
 \]
 Moreover, we assume that the set of functions
 \[
  \mathrm{W}=\{v\in L^2:\mathcal{W}(v)<\infty\}
 \]
 is a Banach space such that
 \[
  C_0^{\infty}\hookrightarrow\mathrm{W}\hookrightarrow L^2\qquad\mbox{(dense embeddings).}
 \]
 Finally, $\mathcal{W}$ is G\^ateaux differentiable on $\mathrm{W}$ and its derivative $\nabla\mathcal{W}:\mathrm{W}\to\mathrm{W}'$ satisfies
 \[
  \|\nabla\mathcal{W}(v)\|_{\mathrm{W}'}\leq C(1+\mathcal{W}(v)^{\theta}\,),\qquad\forall v\in\mathrm{W},
 \]
 for suitable constants $C\geq0$ and $\theta\in(0,1)$.
\end{itemize}

\begin{remark}
 Assumption {\bf (H)}  is typically satisfied by standard functionals like
 \[
  \mathcal{W}(v)=\frac{1}{p}\int_{\mathbb{R}^n}|\nabla^kv|^p\,dx,\qquad p>1,
 \]
 (with possibly lower order terms) where $\mathrm{W}$ is the space of the $L^2$ functions $v$ with $\nabla^kv\in L^p$.
\end{remark}

\begin{theorem}[Serra$\&$Tilli, \cite{serra:3}]
 \label{theorem:serra-tilli-due}
 Given $w_0,\,w_1\in\mathrm{W}$ and $\varepsilon\in(0,1)$, under assumption {\bf (H)} the functional $F_\varepsilon$ defined in \eqref{eq:functional_homogeneous} has a minimizer $w_{\varepsilon}$ in the space $H_{loc}^2([0,\infty);L^2)$ subject to \eqref{eq:initial_conditions}. Moreover:
 \begin{itemize}
  \item[(a)] Estimates. There exists a constant C, independent of $\varepsilon$, such that
  \[
   \begin{array}{c}
    \displaystyle\int_{\tau}^{\tau+T}\mathcal{W}(w_{\varepsilon}(t,\cdot))\,dt\leq CT,\qquad\forall\tau\geq0,\quad\forall T\geq\varepsilon,\\[.5cm]
    \displaystyle\int_{\mathbb{R}^n}|w_{\varepsilon}'(t,x)|^2\,dx\leq C\quad\mbox{and}\quad\int_{\mathbb{R}^n}|w_{\varepsilon}(t,x)|^2\,dx\leq C(1+t^2),\qquad\forall t\geq0,\\[.5cm]
    \displaystyle\|w_{\varepsilon}\|_{L^{\infty}(\mathbb{R}^+;\mathrm{W}')}\leq C.
   \end{array}
  \]
  \item[(b)] Convergence. Every sequence $w_{\varepsilon_i}$ (with $\varepsilon_i\downarrow0$) admits a subsequence which is convergent, in the weak topology of $H^1((0,T);L^2)$ for every $T>0$, to a function $w$ such that
  \[
   w\in H_{loc}^1([0,\infty);L^2),\qquad w'\in L^{\infty}(\mathbb{R}^+;L^2),\qquad w''\in L^{\infty}(\mathbb{R}^+;\mathrm{W}').
  \]
  Moreover, $w$ satisfies the initial conditions \eqref{eq:initial_conditions}.
  \item[(c)] Energy inequality. Letting
  \begin{equation}
   \label{eq:mechanical_energy}
   \mathcal{E}(t):=\frac{1}{2}\int_{\mathbb{R}^n}|w'(t,x)|^2\,dx+\mathcal{W}(w(t,\cdot)),
  \end{equation}
  the function $w(t,x)$ satisfies the energy inequality
  \[
   \mathcal{E}(t)\leq\mathcal{E}(0)=\frac{1}{2}\int_{\mathbb{R}^n}|w_1(x)|^2\,dx+\mathcal{W}(w_0)\qquad\mbox{for a.e. }t>0.
  \]
 \end{itemize}
\end{theorem}

Unfortunately, under these assumptions, it is not known whether $w$ satisfies \eqref{eq:wave_general}. Anyway, Serra$\&$Tilli, still in \cite{serra:3}, provided a \emph{sufficient} condition on $\mathcal{W}$ that allows to obtain \eqref{eq:wave_general}.

\begin{theorem}[Serra$\&$Tilli, \cite{serra:3}]
 \label{theorem:serra-tilli-tre}
 Assume that, for some real number $m>0$,
 \begin{equation}
  \label{eq:further_assumption}
  \mathcal{W}(v)=\frac{1}{2}\|v\|_{\dot{H}^m}^2+\sum_{0\leq k<m}\frac{\lambda_k}{p_k}\int_{\mathbb{R}^n}|\nabla^kv(x)|^{p_k}\,dx\qquad(\lambda_k\geq0,\,p_k>1).
 \end{equation}
 Then, assumption {\bf (H)} is fulfilled if $\mathrm{W}$ is the space of those $v\in H^m$ with $\nabla^kv\in L^{p_k}$ $(0\leq k<m)$ endowed with its natural norm.
 
 Moreover, the limit function $w$ obtained via Theorem \ref{theorem:serra-tilli-due} solves, in the sense of distributions, the hyperbolic equation \eqref{eq:wave_general}.
\end{theorem}

\begin{remark}
 We recall that, as usual, $\|v\|_{\dot{H}^m}$ is the $L^2$ norm of $|\xi|^m\,\widehat{v}(\xi)$, where $\widehat{v}$ is the Fourier transform of $v$. The typical case is when $m$ is an integer, so that $\|v\|_{\dot{H}^m}$ reduces to $\|\nabla^mv\|_{L^2}$.
\end{remark}

\subsection{Examples.}

In addition to \eqref{eq:wave_serra-tilli} there are many other second order hyperbolic equations that can be investigated using the approach suggested by \textsc{Theorem} \ref{theorem:serra-tilli-due}$\&$\textsc{Theorem} \ref{theorem:serra-tilli-tre}. We briefly recall here some of the most significant ones (for a complete discussion we refer the reader to \cite{serra:3}):

\begin{itemize}
 \item[1.] {\it Nonlinear vibrating--beam equation}:
 \[
  w''=-\Delta^2w+\Delta_pw-|w|^{q-2}\,w\qquad(p,\,q>1).
 \]
 Here $\mathcal{W}$ is defined by
\[
 \mathcal{W}(v)=\int_{\mathbb{R}^n}\left(\frac{1}{2}|\Delta v|^2+\frac{1}{p}|\nabla v|^p+\frac{1}{q}|v|^q\right)\,dx
\]
and $\mathrm{W}=\{v\in H^2:\nabla v\in L^p,v\in L^q\}$.

\item[2.] {\it Wave equation with fractional Laplacian}:
\[
 w''=-(-\Delta)^s\qquad(0<s<1).
\]
Here $\mathcal{W}$ is defined by
\[
 \mathcal{W}(v)=c_{n,s}\int_{\mathbb{R}^n\times\mathbb{R}^n}\frac{|v(x)-v(y)|^2}{|x-y|^{n+2s}}\,dx\,dy
\]
(which is, for a proper choice of $c_{n,s}$, the natural energy associated to the fractional Laplacian) and $\mathrm{W}=H^s$.

\item[3.] {\it Sine--Gordon equation}:
\[
 w''=\Delta w-\sin w.
\]
Here $\mathcal{W}$ is defined by
\[
 \mathcal{W}(v)=\int_{\mathbb{R}^n}\left(\frac{1}{2}|\nabla v|^2+1-\cos v\right)\,dx
\]
and $\mathrm{W}=H^1$.

\item[4.] {\it Wave equation with $p$--Laplacian }:
\[
 w''=\Delta_p w.
\]
Here $\mathcal{W}$ is defined by
\[
 \mathcal{W}(v)=\frac{1}{p}\int_{\mathbb{R}^n}|\nabla v|^p\,dx
\]
and $\mathrm{W}=\{v\in L^2:\nabla v\in L^p\}$. 
\end{itemize}

Note that in the cases of the Sine--Gordon and the $p$--Laplacian equation, the functional $\mathcal{W}$ satisfies assumption {\bf (H)}, but not assumption \eqref{eq:further_assumption}. Consequently, one could not apply \textsc{Theorem} \ref{theorem:serra-tilli-tre} to this cases. In the Sine--Gordon case, however, since the functional is quadratic in the higher order space derivatives, one can prove an analogous of \textsc{Theorem} \ref{theorem:serra-tilli-tre}. On the contrary, it is an open problem whether this can be done also for the case of the $p$--Laplacian.

\section{Extension to nonhomogeneous equations.}

The natural further extension is the addition of a general forcing term at the right--hand side of \eqref{eq:wave_general}, that is, the study of the Cauchy problem associated with the nonhomogeneous equation
\begin{equation}
 \label{eq:wave_nonhomogeneous}
  w''(t,x)=-\nabla\mathcal{W}(w(t,\cdot))(x)+f(t,x).
\end{equation}
The proper choice for the functional $F_\varepsilon$ in this case is given by
\begin{equation}
 \label{eq:functional_nonhomogeneous}
 F_\varepsilon(u)=\int_0^te^{-t/\varepsilon}\bigg(\int_{\mathbb{R}^n}\frac{\varepsilon^2|u''(t,x)|}{2}\,dx+\mathcal{W}(u(t,\cdot))-\int_{\mathbb{R}^n}f_\varepsilon(t,x)u(t,x)\,dx\bigg)\,dt,
\end{equation}
where $(f_\varepsilon)$ is a sequence suitably converging to $f$.

This issue has been the topic of the doctoral dissertation of the author and is extensively investigated in \cite{tentarelli}. Here we just announce the result.

\begin{theorem}[Tentarelli$\&$Tilli, \cite{tentarelli}]
 \label{theorem:tentarelli-tilli}
 Let $\mathcal{W}$ be a functional satisfying assumption {\bf (H)} and $w_0,\,w_1\in\mathrm{W}$. Let also $f\in L_{loc}^2([0,\infty),L^2)$. Then, there exists a sequence $(f_\varepsilon)$, converging to $f$ in $L^2([0,T];L^2)$ for all $T>0$, such that:
 \begin{itemize}
  \item[(a)] Minimizers. For every $\varepsilon\in(0,1)$, the functional $F_{\varepsilon}$ defined by \eqref{eq:functional_nonhomogeneous} has a minimizers $w_{\varepsilon}$ in the class of functions in $H_{loc}^2([0,\infty);L^2)$ that are subject to \eqref{eq:initial_conditions}.
  \item[(b)] Estimates. There exist two positive constants $C_t,\,C_{\tau,T}$, depending on $t,\,\tau$ and $T$ (in a continuous way), but independent of $\varepsilon$, such that
  \[
  \int_{\mathbb{R}^n}|w_{\varepsilon}'(t,x)|^2\,dx\leq C_{t},\qquad\int_{\mathbb{R}^n}|w_{\varepsilon}(t,x)|^2\,dx\leq C_{t},\qquad\forall t\geq0,
 \]
 \[
  \int_{\tau}^{\tau+T}\mathcal{W}(w_{\varepsilon}(t,\cdot))\,dt\leq C_{\tau,T},\qquad\forall\tau\geq0,\quad\forall T>\varepsilon,
 \]
 \[
  \int_0^t\|w_{\varepsilon}''(s)\|_{\mathrm{W}'}^2\,ds\leq C_{t},\qquad\forall t\geq0.
 \]
 \item[(c)] Convergence. Every sequence $w_{\varepsilon_i}$ (with $\varepsilon_i\downarrow0$) admits a subsequence which is convergent in the weak topology of $H^1([0,T];L^2)$, for every $T>0$, to a function $w$ that satisfies \eqref{eq:initial_conditions} (where the latter is meant as an equality in $\mathrm{W}'$). In addition,
 \[
  w'\in L_{loc}^{\infty}([0,\infty);L^2)\qquad\text{and}\qquad w''\in L_{loc}^2([0,\infty);\mathrm{W}').
 \]
 \item[(d)] Energy inequality. Letting $\mathcal{E}$ be again the mechanical energy defined by \eqref{eq:mechanical_energy}, there results
 \begin{equation}
  \label{eq:energy_estimate}
  \mathcal{E}(t)\leq\left(\sqrt{\mathcal{E}(0)}+\sqrt{\frac{t}{2}\int_0^t\int_{\mathbb{R}^n}|f(s,x)|^2\,dx\,ds\,}\right)^{2},\qquad\text{for a.e.}\quad t\geq0.
 \end{equation}
 \item[(e)] Solution of \eqref{eq:wave_nonhomogeneous}. Assuming, furthermore, that for some real numbers $m>0,\,\lambda_k\geq0$ and $p_k>1$, $\mathcal{W}$ satisfies \eqref{eq:further_assumption}, then the limit function $w$ solves \eqref{eq:wave_nonhomogeneous}.
 \end{itemize}
\end{theorem}

Some comments are in order. First, we point out that the estimate on the mechanical energy established by \eqref{eq:energy_estimate} is the same that one can find applying a formal {\it Gr\"onwall--type} argument to \eqref{eq:wave_nonhomogeneous}. In addition, setting $f\equiv0$, the results of \textsc{Theorem} \ref{theorem:tentarelli-tilli} recover exactly those of \textsc{Theorem} \ref{theorem:serra-tilli-due}$\&$\textsc{Theorem} \ref{theorem:serra-tilli-tre}, thus showing that our extension to nonhomogeneous equations is consistent.

On the other hand, the Euler--Lagrange equation satisfied by the minimizers of $F_\varepsilon$ (which is analogous to \eqref{eq:euler_lagrange}) suggests to work directly with $f$ in place of $f_{\varepsilon}$ in \eqref{eq:functional_nonhomogeneous}. However, this gives rise to several issues in establishing the requested a priori estimates on $F_{\varepsilon}(w_{\varepsilon})$. On the contrary, a proper choice of $(f_{\varepsilon})$ allows one to adapt the De Giorgi approach under the sole assumption $f\in L_{loc}^2([0,\infty);L^2)$, which is the usual one in the search of solutions of {\it finite energy} for \eqref{eq:wave_nonhomogeneous}. In particular, the detection of a proper (topology and) ``speed of convergence'' for $f_{\varepsilon}$ to $f$ is one of the main issues in the extension to nonhomogeneous problems.

Finally, it is worth to outline briefly the main difference between the homogeneous and the nonhomogeneous case: the estimates on the sequence $(w_{\varepsilon})$ are no longer {\it global} in time. This occurs since the presence of $f$ drops all the uniform bounds deduced in \cite{serra:3} and allows to establish estimates that are either independent of $\varepsilon$ or independent of $t$. In particular, the presence of the forcing term entails that the quantity
\[
 E_{\varepsilon}(t):=\frac{1}{2}\int_{\mathbb{R}^n}|w_{\varepsilon}'(t,x)|^2\,dx+\int_t^{\infty}\varepsilon^{-2}\,e^{-(s-t)/\varepsilon}\,(s-t)\,\mathcal{W}(w_{\varepsilon}(s,\cdot))\,ds
\]
is not decreasing (as in the homogeneous case) and not even uniformly bounded with respect to both $\varepsilon$ and $t$. This function, that we call \emph{approximate energy}, is a {\it formal} approximation of the mechanical energy $\mathcal{E}$ and the investigation of its behavior is the main point of our approach, since it provides the a priori estimates on the minimizers $w_{\varepsilon}$. Consequently, the fact that it admits only estimates on bounded intervals is the reason for which the inequalities in {\it (b)} are no longer global.

This transition ``from global to local'' of the a priori estimates affects the regularity of the limit function $w$, but fortunately does not rule out the possibility of extending the De Giorgi approach. Actually, the proofs of the energy inequality, the initial conditions and \eqref{eq:wave_nonhomogeneous} do not require any global estimate on the sequence of minimizers (even in the homogeneous case).

Moreover, we point out that the choice of the sequence $(f_{\varepsilon})$ is crucial also for establishing the proper estimate on $E_{\varepsilon}(t)$; in particular, for establishing {\it causal} estimates for a quantity which is {\it a--causal} by definition.

\section{A further extension: dissipative equations.}

Finally, it worth recalling that \cite{serra:3} also shows that an approach {\it \`a la De Giorgi} is available also for {\it dissipative} homogeneous wave equations of the type
\[
 w''(t,x)=-\nabla\mathcal{W}(w(t,\cdot))(x)-\nabla\mathcal{G}(w'(t,\cdot))(x),
\]
where $\mathcal{G}$ is a {\it quadratic form} defined on a suitable Hilbert space. A typical example is given by the {\it Telegraph} equation
\[
 w''=\Delta w-|w|^{p-2}\,w-w'\qquad(p>1)
\]
(just setting $\mathcal{W}(v)=\int_{\mathbb{R}^n}\left(\frac{1}{2}|\nabla v|+\frac{1}{p}|v|^p\right)\,dx$ and $\mathcal{G}(v)=\frac{1}{2}\int_{\mathbb{R}^n}|v|^2\,dx$).

As in the non--dissipative case, also here it is natural to wonder if an extension to the nonhomogeneous case, namely
\[
 w''(t,x)=-\nabla\mathcal{W}(w(t,\cdot))(x)-\nabla\mathcal{G}(w'(t,\cdot))(x)+f(t,x),
\]
is possible. The answer is again positive and this issue will be treated in a forthcoming paper by the author, as well.

\end{document}